\documentclass[12pt,reqno,a4paper]{amsart} 

\usepackage{amssymb,latexsym}
\usepackage{cite}
\usepackage[height=200mm,width=140mm]{geometry}
\theoremstyle{plain}
\newtheorem{theorem}{Theorem}
\newtheorem{lemma}{Lemma}

\theoremstyle{definition}

\theoremstyle{remark}

\begin{document}
\title[Finsler-ArealSpaces]{On the structure of Finsler and areal spaces} 

\author{E. Tanaka}
	\address{Erico Tanaka\\ Ochanomizu University\\ Department of Physics\\ 2-1-1 Ootsuka Bunkyo\\112-0086 Tokyo \\Japan}
\email{erico@cosmos.phys.ocha.ac.jp}
	\address{Palacky University\\ Department of Mathematics\\ 17. Listopadu 12\\ 77146 
 Olomouc\\ Czech Republic}

\author{D. Krupka}
	\address{Demeter Krupka\\ Lepage Research Institute\\ Slatinice \\ 78342
Olomouc\\ Czech Republic}
\email{demeter.krupka@lepageri.eu}
	\address{Beijing Institute of Technology\\ Department of Mathematics\\No.5 South Zhongguancun Street\\ 100081
  Beijing\\ China}
\address{The University of Ostrava\\ Department of Mathematics\\30. dubna 22\\ 70103
  Ostrava\\ Czech Republic}

\begin{abstract}
We study underlying geometric structures for integral variational functionals, 
depending on submanifolds of a given manifold. Applications include (first order) 
variational functionals of Finsler and areal geometries with integrand the Hilbert 1-form, 
and admit immediate extensions to higher-order functionals. 
\end{abstract}

\subjclass{49Q99, 53C60, 58E30}

\keywords{Submanifold, $k$-vector, Grassmann fibration, Hilbert form}

\maketitle

\section{Introduction}
This paper is a contribution to the theory of integral variational functionals, 
depending on submanifolds of a given manifold $X$. The theory is based on geometric notions, 
such as the bundles of (skew-symmetric) multivectors, and Grassmann fibrations. 
Conceptually, it extends local parametric integrals of Finsler-Kawaguchi and areal geometries 
(see e.g. Chern, Chen, Lam~\cite{ChernChenLam}, Davies~\cite{Davies}, Kawaguchi~\cite{AK5}, 
and Tamassy~\cite{Tamassy4}) 
to global functionals, depending on (global) submanifolds. 
In Section 2 we summarize integration theory of differential forms along submanifolds. 
Section 3 is devoted to vector bundles of $k$-vectors; we show how mappings of Euclidean spaces into manifolds 
({\it parametrisations}) can be lifted to the bundles of $k$-vectors. 
In Section 4 we introduce, using the Pl\"ucker embedding, 
underlying spaces for parameter-invariant variational problems, the Grassmann fibrations. 
In Section 5 we show that any $k$-form on the Grassmann fibration defines an integral variational functional, 
depending on $k$-dimensional submanifolds. An example is the {\it Hilbert form}, a well-known first-order construction 
in Finsler geometry and its generalisations (Chern, Chen, Lam \cite{ChernChenLam}, 
Crampin, Saunders~\cite{CrampinSaunders}). 

It should be pointed out that the theory can be further generalised. 
To this purpose one should consider higher-order Grassmann fibrations, 
endowed with Lagrangians, satisfying relevant homogeneity conditions 
(Zermelo conditions, see e.g. Saunders~\cite{Saunders3}, and Urban and Krupka~\cite{Krupka-Urban}). 

\section{Integration over submanifolds}
Let $X$ be an $n$-dimensional manifold, $S$ a  subset of $X$, ${x_0} \in S$ a point. 
A chart $(U,\varphi )$,$\varphi  = ({x^i})$, at ${x_0}$ is a {\it submanifold chart} for $S$, 
if there exists a non-negative integer $k \leqslant n$ such that 
$\varphi (U \cap S) = \{ x \in U|{x^{k + 1}}(x) = {c_1},{x^{k + 2}}(x) = {c_2}, \ldots ,{x^n}(x) = {c_{n - k}}\} $. 
If such a chart exists, we say that $S$ is a {\it submanifold} of $X$ at the point ${x_0}$; 
$k$ is the {\it dimension} of $S$ at ${x_0}$. If such a submanifold chart exists at every point of $X$, 
we say $S$ is a {\it submanifold} of $X$ and call $k$ the {\it dimension} of $S$. 

Denote by $({t^1},{t^2}, \ldots ,{t^n})$ the canonical coordinates on the Euclidean space 
${R^n}$, and $R_{( - )}^n = \{ {t_0} \in {R^n}|{t^n}({t_0}) \leqslant 0\} $, $\partial R_{( - )}^n 
= \{ {t_0} \in R_{( - )}^n|{t^n}({t_0}) = 0\} $. $R_{( - )}^n$ is the {\it halfspace} of ${R^n}$, 
$\partial R_{( - )}^n$ is the {\it boundary} of $R_{( - )}^n$. Let $\Omega $ be a non-void subset of $X$, 
and ${x_0} \in \Omega $ a point. A chart $(U,\varphi )$ at ${x_0}$ is said to be {\it adapted} to $\Omega $, 
if the set $\varphi (U \cap \Omega )$ is an open set in $R_{( - )}^n$. $\Omega $ is a {\it piece} of $X$, 
if it is compact and each point $x \in \Omega $ admits a chart, adapted to $\Omega $. 

Let $\eta $ be a $k$-form on $X$. Our aim now will be to introduce an integral of $\eta $ on a piece of a 
$k$-dimensional submanifold $S$ ({\it $k$-piece} of a $X$). Express $\eta $ in a submanifold chart 
$(U,\varphi )$, $\varphi  = ({x^i})$, as 
$\eta  = {\eta _{{i_1}{i_2} \ldots {i_k}}}d{x^{{i_1}}} 
\wedge d{x^{{i_2}}} \wedge  \ldots  \wedge d{x^{{i_k}}}.$
Then restricting $\eta $ to $S$ 
we get from the equations $x^{k + 1} = 0$, ${x^{k + 2}} = 0$, $\cdots$ , ${x^n} = 0$ 
\begin{align}
\eta  = fd{x^1} \wedge d{x^2} \wedge  \ldots  \wedge d{x^k},
\end{align}
where we write $f = f({x^{{i_1}}},{x^{{i_2}}}, \ldots ,{x^{{i_k}}})$ 
for the component of $\eta $ restricted to $S$. 
From now on we suppose that $S$ is {\it orientable}, and is endowed with an orientation 
${\operatorname{Or} _S}X$; only submanifold charts on $X$ belonging to ${\operatorname{Or} _S}X$ are used. 
The integral of $\eta $ on a compact set $\Omega  \subset S$ is defined in a standard way. 
There exist a finite family $\{ ({U_1},{\varphi _1}),({U_2},{\varphi _2}), \ldots ,({U_N},{\varphi _N})\} $ 
of submanifold charts on $X$, such that the family $\{ {U_1} \cap S,{U_2} \cap S, \ldots ,{U_N} \cap S\} $ 
covers $\Omega $. Let $\{ {\chi _1},{\chi _2}, \ldots ,{\chi _N}\} $ be a partition of unity, 
subordinate to this covering. Then 
\begin{align}
\int_\Omega  \eta   = \sum\limits_{j = 1}^N {\int_{\operatorname{supp} {\chi _j} \cap \Omega } {{\chi _j}\eta } } .
\end{align}
The following basic properties of the integral are needed in the calculus of variations. 
\begin{lemma} {\bf (Transformation of integration domain)} 
Let $X$ and $Y$ be two smooth n-dimensional oriented manifolds, 
$\alpha :X \to Y$ an orientation-preserving diffeomorphism. 
Then for any compact set $\Omega  \subset S$ and any continuous differential $n$-form on $Y$
\begin{align}
\int_\Omega  \eta   = \int_{{\alpha ^{ - 1}}(\Omega )} {\alpha *\eta } .
\end{align}
\end{lemma}
\begin{lemma} {\bf (Leibniz rule)} 
Let $X$ be an oriented n-dimensional manifold, ${\eta _t}$ a family of $n$-forms on $X$, 
differentiable on a real parameter $t$, $\Omega  \subset S$ a compact set. 
Then the function 
\begin{align}
I \ni t \to \int_\Omega  {{\eta _t}}  \in R
\end{align}
is differentiable, and 
\begin{align}
\frac{d}{{dt}}\int_\Omega  {{\eta _t}}  = \int_\Omega  {\frac{{d{\eta _t}}}{{dt}}} .
\end{align}
\end{lemma}
\begin{lemma} {\bf (The Stokes' formula) }  
Let $X$ be an $n$-dimensional manifold, $S$ a $k$-dimensional oriented submanifold of $X$, 
$\eta $ a $(k - 1)$-form on $X$. 
Let $\Omega $ be a piece of $S$ with boundary $\partial \Omega $, endowed with induced orientation. Then 
\begin{align}
\int_{\partial \Omega } \eta   = \int_\Omega  {d\eta } .
\end{align}
\end{lemma}

\section{Bundles of $k$-vectors}
Let $X$ be an $n$-dimensional manifold, 
${\Lambda ^k}{T_x}X$ the $k$-exterior product of the tangent space 
${T_x}X$, $x \in X$ a point. We put 
\begin{align}
{\Lambda ^k}TX = \bigcup\limits_{x \in X} {{\Lambda ^k}{T_x}X} .
\end{align}
This set has a natural vector bundle structure over $X$, with type fibre ${\Lambda ^k}{R^n}$. 
We denote by ${\tau ^k}$ the vector bundle projection of ${\Lambda ^k}TX$.
 
Let $X$ (resp. $Y$) be a smooth manifold of dimension $n$ (resp. $m$), 
and let $f:X \to Y$ be a differentiable mapping. 
Choose a point $x \in X$ and a $k$-vector $\Xi  \in {\Lambda ^k}{T_x}X$. 
Then choose a chart $(U,\varphi )$, $\varphi  = ({x^i})$, at $x$ and a chart $(V,\psi )$, 
$\psi  = ({y^\sigma })$, at $f(x) \in Y$ such that $f(U) \subset V$. Expressing $\Xi $ 
in components and setting 
\begin{align}
&{\Lambda ^k}{T_x}f \cdot \Xi  = \frac{1}{{{{(k!)}^2}}}{\left( {\frac{{\partial {y^{{\sigma _1}}}
f{\varphi ^{ - 1}}}}{{\partial {x^{{i_1}}}}}} \right)_{\varphi (x)}}{\left( {\frac{{\partial {y^{{\sigma _2}}}
f{\varphi ^{ - 1}}}}{{\partial {x^{{i_2}}}}}} \right)_{\varphi (x)}} \ldots {\left( 
{\frac{{\partial {y^{{\sigma _k}}}f{\varphi ^{ - 1}}}}{{\partial {x^{{i_k}}}}}} \right)_{\varphi (x)}} \hfill \nonumber \\
&\quad  \cdot {\Xi ^{{i_1}{i_2} \ldots {i_k}}}{\left( {\frac{\partial }{{\partial {y^{{\sigma _1}}}}}} 
\right)_{f(x)}} \wedge {\left( {\frac{\partial }{{\partial {y^{{\sigma _2}}}}}} \right)_{f(x)}} \wedge 
 \ldots  \wedge {\left( {\frac{\partial }{{\partial {y^{{\sigma _k}}}}}} \right)_{f(x)}}, 
\end{align}
we get a $k$-vector ${\Lambda ^k}{T_x}f \cdot \Xi  \in {\Lambda ^k}{T_{f(x)}}Y$, 
and a vector bundle homomorphism ${\Lambda ^k}Tf:{\Lambda ^k}TX \to {\Lambda ^k}TY$ over  $f$ (the {\it lift} of $f$). 

It is easily seen that differentiable mappings of a Euclidean space into a manifold can be canonically lifted 
to the bundles of $k$-vectors. To this purpose we use the {\it canonical $k$-vector field} on ${R^n}$
\begin{align}
{R^n} \ni t \to \theta (t) = \frac{1}{{k!}}{\varepsilon ^{{i_1}{i_2} \ldots {i_k}}}
{\left( {\frac{\partial }{{\partial {t^{{i_1}}}}}} \right)_t} \wedge 
{\left( {\frac{\partial }{{\partial {t^{{i_2}}}}}} \right)_t} \wedge  
\ldots  \wedge {\left( {\frac{\partial }{{\partial {t^{{i_k}}}}}} \right)_t} \in {\Lambda ^k}T{R^n}.
\end{align}
Identifying ${\Lambda ^k}T{R^n}$ with ${R^n} \times {\Lambda ^k}{R^n}$, 
the canonical section becomes the mapping $t \to (t,{\varepsilon ^{{i_1}{i_2} \ldots {i_k}}})$. 

Consider a differentiable mapping $f:U \to Y$, where U is an open subset of ${R^n}$. 
For any point $t \in U$, ${\Lambda ^k}{T_t}f \cdot \theta (t)$ is an element of the vector space 
${\Lambda ^k}{T_{f(t)}}Y$. We get the {\it canonical lift} ${\Lambda ^k}f$ of $f$ to ${\Lambda ^k}TY$, defined by 
\begin{align}
{\Lambda ^k}f = {\Lambda ^k}Tf \cdot \theta .
\end{align}
The canonical lift of the {\it parametrisation} $U \ni t \to ({\psi ^{ - 1}} \circ {\iota _{k,m}})(t) \in V \cap S$ 
is expressed in a chart $(V,\psi )$, $\psi  = ({y^\sigma })$, as
\begin{align}
&{\Lambda ^k}({\psi ^{ - 1}} \circ {\iota _{k,m}})(t) \hfill \nonumber \\
&\quad  = {\left( {\frac{\partial }{{\partial {y^1}}}} \right)_{{\psi ^{ - 1}} 
\circ {\iota _{k,m}}(t)}} \wedge {\left( {\frac{\partial }{{\partial {y^2}}}} \right)_{{\psi ^{ - 1}} 
\circ {\iota _{k,m}}(t)}} \wedge  \ldots  \wedge {\left( {\frac{\partial }{{\partial {y^k}}}} 
\right)_{{\psi ^{ - 1}} \circ {\iota _{k,m}}(t)}}. \label{form5}
\end{align}
Formula (\ref{form5}) also defines the mapping 
$V \ni y \to ({\Lambda ^k}\psi )(y) \in {({\tau ^k})^{ - 1}}(V)$ by 
\begin{align}
{\Lambda ^k}\psi  = {\Lambda ^k}({\psi ^{ - 1}} \circ {\iota _{k,m}}) \circ {\operatorname{pr} _{m,k}}\psi ,
\end{align}
the {\it canonical section along} $S$, associated with $(V,\psi )$. ${\Lambda ^k}\psi $ is expressed by 
\begin{align}
&({y^1},{y^2}, \ldots ,{y^k},{y^{k + 1}},{y^{k + 2}}, \ldots ,{y^m}) \to {\Lambda ^k}({\psi ^{ - 1}} 
\circ {\iota _{k,m}})({y^1},{y^2}, \ldots ,{y^k}) \hfill \nonumber \\
&\quad  = (({y^1},{y^2}, \ldots ,{y^k},0,0, \ldots ,0),(1,0,0, \ldots ,0)). 
\end{align}
Writing in the multi-index notation $({({\tau ^r})^{ - 1}}(V),\Phi )$, 
$\Phi  = ({\dot y^I})$, and setting ${I_0} = (12 \ldots k)$, 
we get the image of this mapping as a subset of 
${({\tau ^r})^{ - 1}}(V)$, defined by the equations 
${y^{k + 1}} = 0$, ${y^{k + 2}} = 0$, $ \ldots $, ${y^m} = 0$, 
${\dot y^{{I_0}}} = 1$, ${\dot y^I} = 0$, $I \ne {I_0}$. 
\begin{lemma}
Let $(V,\psi )$, $\psi  = ({y^\sigma })$, and $(\bar V,\bar \psi )$, 
$\bar \psi  = ({\bar y^\sigma })$, be two charts on $Y$, adapted to $S$, such that $V \cap \bar V \ne \O $. 
\begin{enumerate}
\item The canonical sections along $S$ satisfy 
\begin{align}
{\Lambda ^k}\bar \psi  = \det {\left( {\frac{{\partial {y^i}}}{{\partial {{\bar y}^j}}}} 
\right)_{\psi (y)}}_{\bar \psi (y)}{\Lambda ^k}\psi .
\end{align}

\item The differential forms ${d y^\sigma }$ and $d{y^I}$ satisfy 
$({\Lambda ^k}\psi )^* d{y^i} = d{y^i}$, $1 \leqslant i \leqslant k$, 
$({\Lambda ^k}\psi )^* d{y^\nu } = 0$, $k + 1 \leqslant \nu  \leqslant m$, 
$({\Lambda ^k}\psi )^* d{\dot y^I} = 0$. 
In particular, on the set $V \cap \bar V$,
\begin{align}
& ({\Lambda ^k}\bar \psi )^*d{{\bar y}^1} \wedge d{{\bar y}^2} \wedge  \ldots  \wedge d{{\bar y}^k} \hfill \nonumber \\
&  \quad  = \det {\left( {\frac{{\partial {{\bar y}^i}}}{{\partial {y^j}}}} \right)_{\psi (y)}}({\Lambda ^k}\psi )*d{y^1} \wedge d{y^2} \wedge  \ldots  \wedge d{y^k}. 
\end{align}
\end{enumerate}
\end{lemma}

\section{Grassmann fibrations}
Consider the vector bundle ${\Lambda ^k}TY$ and the subset $\Lambda _0^kTY \subset {\Lambda ^k}TY$, 
consisted of {\it non-zero} $k$-vectors. We have an equivalence relation on $\Lambda _0^kTY$ 
``{\it ${\Xi _1}$ is equivalent with ${\Xi _2}$, if there exists a real number $\lambda  > 0$ 
such that ${\Xi _1} = \lambda {\Xi _2}$}''. 
The quotient set has the structure of a fibration over $Y$, called the {\it Grassmann fibration} of degree $k$, 
and denoted by ${G^k}Y$. 
	
To describe the structure of the set ${G^k}Y$, we proceed similarly as in the case of classical projective spaces. 
If in a chart $(V,\psi )$, $\psi  = ({y^\sigma })$, 
\begin{align}
{\Xi _i} = \frac{1}{{k!}}\Xi _i^{{\sigma _1}{\sigma _2} \ldots {\sigma _k}}
{\left( {\frac{\partial }{{\partial {y^{{\sigma _1}}}}}} \right)_y} \wedge 
{\left( {\frac{\partial }{{\partial {y^{{\sigma _2}}}}}} \right)_y} \wedge  
\ldots  \wedge {\left( {\frac{\partial }{{\partial {y^{{\sigma _k}}}}}} \right)_y},\quad i = 1,2,
\end{align}
are two nonzero $k$-vectors, then ${\Xi _1}$ is equivalent with ${\Xi _2}$ if and only if in this chart, 
$\Xi _1^{{\sigma _1}{\sigma _2} \ldots {\sigma _k}} = \lambda \Xi _2^{{\sigma _1}{\sigma _2} 
\ldots {\sigma _k}}$ for some $\lambda  > 0$ and all 
${\sigma _1},{\sigma _2}, \ldots ,{\sigma _k}$. 
We denote ${V^{{\nu _1}{\nu _2} \ldots {\nu _k}}} = \{ \Xi  
\in {({\tau ^k})^{ - 1}}(V) \,|\, {\Xi ^{{\nu _1}{\nu _2} \ldots {\nu _k}}} > 0\} $. 
Then a $k$-vector belonging to the set ${V^{{\nu _1}{\nu _2} \ldots {\nu _k}}} \subset \Lambda _0^k TY$ 
can be expressed by 
\begin{align}
&\Xi  = {\Xi ^{{\nu _1}{\nu _2} \ldots {\nu _k}}}{\left( {\frac{\partial }{{\partial {y^{{\nu _1}}}}}} 
\right)_y} \wedge {\left( {\frac{\partial }{{\partial {y^{{\nu _2}}}}}} \right)_y} \wedge  \ldots  \wedge 
{\left( {\frac{\partial }{{\partial {y^{{\nu _k}}}}}} \right)_y} \hfill \nonumber \\
& \quad  + \frac{1}{{k!}}\sum\limits_{({\tau _1}{\tau _2} \ldots {\tau _k}) 
\ne ({\nu _1}{\nu _2} \ldots {\nu _k})} {{\Xi ^{{\tau _1}{\tau _2} \ldots {\tau _k}}}} 
{\left( {\frac{\partial }{{\partial {y^{{\tau _1}}}}}} \right)_y} \wedge 
{\left( {\frac{\partial }{{\partial {y^{{\tau _2}}}}}} \right)_y} \wedge  
\ldots  \wedge {\left( {\frac{\partial }{{\partial {y^{{\tau _k}}}}}} \right)_y}
\end{align}
(no summation through ${\nu _1},{\nu _2}, \ldots ,{\nu _k}$). 
Denoting by $\operatorname{sgn} {\Xi ^{{\nu _1}{\nu _2} \ldots {\nu _k}}}$ 
the sign of the component ${\Xi ^{{\nu _1}{\nu _2} \ldots {\nu _k}}}$, we can write 
${\Xi ^{{\nu _1}{\nu _2} \ldots {\nu _k}}} = \operatorname{sgn} {\Xi ^{{\nu _1}{\nu _2} \ldots {\nu _k}}} \cdot | 
{\Xi ^{{\nu _1}{\nu _2} \ldots {\nu _k}}}|$ and 
\begin{align}
& \Xi  = \operatorname{sgn} {\Xi ^{{\nu _1}{\nu _2} \ldots {\nu _k}}} \cdot |{\Xi ^{{\nu _1}{\nu _2} 
\ldots {\nu _k}}}|{\left( {\frac{\partial }{{\partial {y^{{\nu _1}}}}}} \right)_y} 
\wedge {\left( {\frac{\partial }{{\partial {y^{{\nu _2}}}}}} \right)_y} \wedge  \ldots  \wedge 
{\left( {\frac{\partial }{{\partial {y^{{\nu _k}}}}}} \right)_y} \hfill \nonumber \\
& \quad  + \frac{{|{\Xi ^{{\nu _1}{\nu _2} \ldots {\nu _k}}}|}}{{k!}}{\sum {\frac{{{\Xi ^{{\tau _1}{\tau _2} 
\ldots {\tau _k}}}}}{{|{\Xi ^{{\nu _1}{\nu _2} \ldots {\nu _k}}}|}}{{\left( 
{\frac{\partial }{{\partial {y^{{\tau _1}}}}}} \right)}_y} 
\wedge {{\left( {\frac{\partial }{{\partial {y^{{\tau _2}}}}}} \right)}_y} \wedge  \ldots  \wedge 
\left( {\frac{\partial }{{\partial {y^{{\tau _k}}}}}} \right)} _y},
\end{align}
with the summation through 
$({\tau _1}{\tau _2} \ldots {\tau _k}) \ne ({\nu _1}{\nu _2} \ldots {\nu _k})$. But 
$\operatorname{sgn} {\Xi ^{{\nu _1}{\nu _2} \ldots {\nu _k}}} = 1$, 
so we see the class of $\Xi $ can be represented as   
\begin{align}
& [\Xi ] = {\left( {\frac{\partial }{{\partial {y^{{\nu _1}}}}}} \right)_y} 
\wedge {\left( {\frac{\partial }{{\partial {y^{{\nu _2}}}}}} \right)_y} \wedge  
\ldots  \wedge {\left( {\frac{\partial }{{\partial {y^{{\nu _k}}}}}} \right)_y} \hfill \nonumber \\
&  \quad  + \frac{1}{{k!}}{\sum {\frac{{{\Xi ^{{\tau _1}{\tau _2} \ldots {\tau _k}}}}}{{{\Xi ^{{\nu _1}{\nu _2} 
\ldots {\nu _k}}}}}{{\left( {\frac{\partial }{{\partial {y^{{\tau _1}}}}}} \right)}_y} 
\wedge {{\left( {\frac{\partial }{{\partial {y^{{\tau _2}}}}}} \right)}_y} \wedge  \ldots  \wedge 
\left( {\frac{\partial }{{\partial {y^{{\tau _k}}}}}} \right)} _y}.
\end{align}
We set for any $\Xi  \in {V^{{\nu _1}{\nu _2} \ldots {\nu _k}}}$
\begin{align}
&  {w^\sigma }(\Xi ) = {y^\sigma }(\Xi ),\quad {w^{{\nu _1}{\nu _2} \ldots {\nu _k}}}(\Xi ) 
= {{\dot y}^{{\nu _1}{\nu _2} \ldots {\nu _k}}}(\Xi ), \hfill \nonumber \\
&  {w^{{\sigma _1}{\sigma _2} \ldots {\sigma _k}}}(\Xi ) = \frac{{{{\dot y}^{{\sigma _1}{\sigma _2} 
\ldots {\sigma _k}}}(\Xi )}}{{{{\dot y}^{{\nu _1}{\nu _2} \ldots {\nu _k}}}(\Xi )}},\quad ({\sigma _1}{\sigma _2} 
\ldots {\sigma _k}) \ne ({\nu _1}{\nu _2} \ldots {\nu _k}). \label{form5-2}
\end{align}
Then the pair $({V^{{\nu _1}{\nu _2} \ldots {\nu _k}}}, 
{\Psi ^{{\nu _1}{\nu _2} \ldots {\nu _k}}})$, 
${\Psi ^{{\nu _1}{\nu _2} \ldots {\nu _k}}} = ({w^\sigma }, {w^{{\nu _1}{\nu _2} \ldots {\nu _k}}}, {w^{{\sigma _1}{\sigma _2} 
\ldots {\sigma _k}}})$, where the indices satisfy 
$({\sigma _1}{\sigma _2} \ldots {\sigma _k}) \ne ({\nu _1}{\nu _2} \ldots {\nu _k})$, 
is a chart on $\Lambda _0^kTY$; we call this chart $({\nu _1}{\nu _2} \ldots {\nu _k})$-{\it associated} 
with $(V,\psi )$. The pair $({V^{{\nu _1}{\nu _2} \ldots {\nu _k}}}, {W^{{\nu _1}{\nu _2} \ldots {\nu _k}}})$, 
${W^{{\nu _1}{\nu _2} \ldots {\nu _k}}} = ({w^\sigma },{w^{{\sigma _1}{\sigma _2} \ldots {\sigma _k}}})$, 
$({\sigma _1}{\sigma _2} \ldots {\sigma _k}) \ne ({\nu _1}{\nu _2} \ldots {\nu _k})$, 
is a fibred chart on ${G^k}Y$. 
Writting formulas (\ref{form5-2}) in a different way, we have the transformation equations 
\begin{align}
{w^\sigma } = {y^\sigma },\quad {w^{{\nu _1}{\nu _2} \ldots {\nu _k}}} = {\dot y^{{\nu _1}{\nu _2} 
\ldots {\nu _k}}},\quad {w^{{\sigma _1}{\sigma _2} \ldots {\sigma _k}}} 
= \frac{{{{\dot y}^{{\sigma _1}{\sigma _2} \ldots {\sigma _k}}}}}{{{{\dot y}^{{\nu _1}{\nu _2} 
\ldots {\nu _k}}}}}.
\end{align}
The projection ${\kappa ^k}:{\Lambda ^k}TY \to {G^k}Y$ of ${\Lambda ^k}TY$ onto 
$G^k Y$ is the Cartesian projection 
$({w^\sigma },{w^{{\nu _1}{\nu _2} \ldots {\nu _k}}}, {w^{{\sigma _1}{\sigma _2} \ldots {\sigma _k}}}) 
\to ({w^\sigma },{w^{{\sigma _1}{\sigma _2} \ldots {\sigma _k}}})$. 
Combining ${\Lambda ^k}({\psi ^{ - 1}}{\iota _{k,m}})$ and 
${\kappa ^k}$ we get the canonical lift of ${\psi ^{ - 1}}{\iota _{k,m}}$ to the Grassmann fibration, 
\begin{align}
{G^k}({\psi ^{ - 1}}{\iota _{k,m}}) = {\kappa ^k} \circ {\Lambda ^k}({\psi ^{ - 1}}{\iota _{k,m}}).
\end{align}
\begin{lemma}
Let $(V,\psi )$, $\psi  = ({y^\sigma })$, and $(\bar V,\bar \psi )$, $\bar \psi  = ({\bar y^\sigma })$, 
be two rectangle charts, adapted to $S$ at a point $y \in Y$. 
Suppose that $(V,\psi )$ and $(\bar V,\bar \psi )$ are consistently oriented. Then 
\begin{align}
{G_k}({\bar \psi ^{ - 1}}{\iota _{k,m}}) = {G_k}({\psi ^{ - 1}}{\iota _{k,m}}).
\end{align}
\end{lemma}
We set
\begin{align}
{G^k}S = \{ [\Xi ] \in {G^k}Y|[\Xi ] = {G^k}({\psi ^{ - 1}}{\iota _{k,m}})({\operatorname{pr} _{m,k}}\psi (y)),y \in S\} .
\end{align}
We associate to a given chart $(V,\psi )$, $\psi  = ({y^\sigma })$, the induced chart 
$({({\tau ^k})^{ - 1}}(V),\Phi )$, $\Phi  = ({y^\sigma },{\dot y^{{\sigma _1}{\sigma _2} \ldots {\sigma _k}}})$, 
on ${\Lambda ^k}TY$; the associated charts on the Grassmann fibration ${G^k}Y$ are $(V_0^{{\nu _1}{\nu _2} 
\ldots {\nu _k}},{W^{{\nu _1}{\nu _2} \ldots {\nu _k}}})$, 
${W^{{\nu _1}{\nu _2} \ldots {\nu _k}}} = ({w^\sigma }, {w^{{\sigma _1}{\sigma _2} \ldots {\sigma _k}}})$, 
$({\sigma _1}{\sigma _2} \ldots {\sigma _k}) \ne ({\nu _1}{\nu _2} \ldots {\nu _k})$. 
Then it is easily seen that each of the charts $(V_0^{{\nu _1}{\nu _2} \ldots {\nu _k}}, {W^{{\nu _1}{\nu _2} 
\ldots {\nu _k}}})$ is adapted to the submanifold ${G^k}S$.
\begin{theorem} \label{th1}
Suppose $S$ is oriented. Then the subset ${G^k}S$ of the Grassmann fibration ${G^k}Y$ is a $k$-dimensional 
oriented submanifold, diffeomorphic with $S$. 
\end{theorem}

Theorem \ref{th1} allows us to integrate over $k$-dimensional submanifolds of $Y$ directly on the 
Grassmann fibration ${G^k}Y$. 

\section{Variational functionals depending on submanifolds}
As before, we denote by ${G^k}S$ (resp. ${G^k}\Omega $) the canonical lift of a $k$-dimensional submanifold 
$S \subset Y$ (resp. $k$-piece $\Omega  \subset Y$) to the Grassmann fibration ${G^k}Y$. 
Denote by ${\Gamma ^k}Y$ the set of all $k$-pieces $\Omega $ of the manifold $Y$. 

Let $\eta $ be a $k$-form on the Grassmann fibration ${G^k}Y$. 
The form $\eta $ defines a {\it variational functional} 
\begin{align}
{\Gamma ^k}Y \ni \Omega  \to {\eta _\Omega }(S) = \int_{{G^k}\Omega } \eta   \in R. \label{vari}
\end{align}
We roughly describe in this paper this construction for $k = 1$, 
representing variational functionals of {\it Finsler geometry} in terms of differential forms 
(cf. Urban and Krupka~\cite{Krupka-Urban2}). 
Consider the tangent bundle ${\Lambda ^1}TY = TY$, a chart $(V,\psi)$, 
$\psi  = ({y^\sigma })$, on $Y$, and the associated chart ${({\tau ^1})^{ - 1}}(V)$, 
$\Psi  = ({y^\sigma },{\dot y^\sigma })$, on $TY$. A function $F:TY \to R$ 
is said to be a {\it Finsler fundamental function}, if it satisfies 
\begin{align}
F(\lambda \xi ) = \lambda F(\xi )
\end{align}
for all tangent vectors $\xi$ and every positive 
$\lambda  \in R$ ({\it homogeneity condition}). The same can be stated in coordinates, requiring that
\begin{align}
F({y^\nu },\lambda {\dot y^\nu }) = \lambda F({y^\nu },{\dot y^\nu }).
\end{align}
\begin{theorem} \label{th2}
\begin{enumerate}
\item For any function $F:TY \to R$, the chart expressions
\begin{align}
\eta  = \frac{{\partial F}}{{\partial {{\dot y}^\nu }}}d{y^\nu }
\end{align}
define a global $1$-form on $TY$. 
\item If $F$ satisfies the homogeneity condition, then $\eta $ is projectable on the Grassmann fibration 
${G^1}TY$. 
\item If $F$ satisfies the homogeneity condition, then for any curve $\zeta :I \to Y$
\begin{align}
({\Lambda ^1}\zeta )*\eta  = (F \circ {\Lambda ^1}\zeta )dt. \label{form5-3}
\end{align}
\end{enumerate}
\end{theorem}
The form $\eta $ (\ref{form5-3}) is known as the {\it Hilbert form} (Chern, Chen and Lam~\cite{ChernChenLam}, 
Crampin and Saunders~\cite{CrampinSaunders}). 
Theorem \ref{th2} (2) characterizes its basic property when F is positive homogeneous:
 namely, in this case the Hilbert form is defined on the {\it Grassmann fibration} ${G^1}TY$. 
 One can also easily verify that $\eta$ is a special case of the {\it Lepage-Cartan form}. 
 This fact completely determines the behaviour of the variational functional (\ref{vari}) 
 under variations of submanifolds, extremal submanifolds, and their invariance properties.

\section*{Acknowledgement} 
Erico Tanaka acknowledges the support by the JSPS Institutional Program for Young Researcher Overseas Visits, 
Palacky University (PrF-2011-022) and Yukawa Institute Computer Facility. 
The second author acknowledges support of National Science Foundation of China (grant No. 109320020). 
He is thankful to the School of Mathematics, Beijing Institute of Technology, 
for kind hospitality and collaboration during his stay in China. 
He also acknowledges support from grant 201/09/0981 of the Czech Science Foundation and from the 
IRSES project GEOMECH (project no. 246981) within the 7th European Community Framework Programme.


\begin{thebibliography}{1}
\expandafter\ifx\csname url\endcsname\relax
  \def\url#1{\texttt{#1}}\fi
\expandafter\ifx\csname urlprefix\endcsname\relax\def\urlprefix{URL }\fi

\bibitem{ChernChenLam}
\textsc{Chern, S.~S., Chen, W.~H., and Lam, K.~S.}: \emph{Lectures on
  Differential Geometry}, World Scientific, 2000.

\bibitem{CrampinSaunders}
\textsc{Crampin, M. and Saunders, D.}: \emph{The hilbert-caratheodory forms for
  parametric multiple integral problems in the calculus of variations}, Acta
  Appl. Math., \textbf{76} (2003), 37--55.

\bibitem{Davies}
\textsc{Davies, E.}: \emph{Areal spaces}, Annali di Matematica Pura ed
  Applicata, \textbf{55} (1961), 63--76.

\bibitem{AK5}
\textsc{Kawaguchi, A.}: \emph{On the concepts and theories of higher order
  spaces}, Periodica Mathematica Hungarica, \textbf{7} (1976), No. 3-4,
  291--299.

\bibitem{Saunders3}
\textsc{Saunders, D.}: \emph{Homogeneous variational problems: a minicourse},
  Communications in Mathematics, \textbf{19} (2010), 91--128.

\bibitem{Tamassy4}
\textsc{Tamassy, L.}: \emph{A contribution to the connection theory of areal
  spaces}, Proc. XXVI Sympos. on Mathematical Physics, Rep. Math. Phys.,
  \textbf{36} (1995), 453--459.

\bibitem{Krupka-Urban2}
\textsc{Urban, Z. and Krupka, D.}: \emph{Variational sequences on fibred
  velocity spaces}.

\bibitem{Krupka-Urban}
\textsc{Urban, Z. and Krupka, D.}: \emph{The zermelo conditions and higher
  order homogeneous functions}, Publ. Math. Debrecen, \textbf{82} (2013),
  59--76.

\end{thebibliography}
\end{document}